\documentclass[11pt,epsf,twoside]{article}

\usepackage[hang,small]{caption2}
\usepackage{graphicx}
\usepackage{rotating}
\usepackage{amssymb,amsmath}
\usepackage[normal]{subfigure}
\usepackage{epsfig}

\usepackage{fancyhdr}
\newcounter{quest}
\newtheorem{question}[quest]{Question}

\def \CC {{\mathbb{C}}}

\def \HH {{\mathbb{H}}}

\begin{document}

\title{\textbf{\huge Playing with the critical point}\\\small An experiment with the Mandelbrot set connectivity}

\author{Alessandro Rosa}
\maketitle

\begin{abstract}
By means of a graphical journey across the Mandelbrot set for the
classic quadratic iterator $f(z):z^2+q$, we illustrate how
connectivity breaks as the seed $z_0$ is no longer at the critical
point of $f(z)$. Finally we suggest an attack to the MLC
conjecture.
\end{abstract}

\section{Introduction}
The Mandelbrot set $\mathcal{M}$ generation relates, in the
original experiment, to a given non-linear and one-parameter map
$f(z,q)$ where both $z,q\in\CC$ and $q$ is the complex parameter
ranging over the whole $\CC$. The key concept of this experiment
refers back to early Fatou's and Julia's ideas on testing Julia
set connectivity under perturbations of the parameter $q$, as well
as on their renown theorem of existence of one critical point $c$
in any basin of attraction and the concept of `post-critical'
set\footnote{Let $\mathcal{C}(f)$ be the set of critical points
for $f$. It is the closure of all forward orbits $f_n(c), n\geq
1$, where $c\in\mathcal{C}(f)$.}, revealing so helpful today to
explore several properties in the dynamics of complex maps. Let
$z_0$ be the seed\footnote{Under this term, one indicates the
first point of any orbit.} at the critical point of $f$. Let $q\in
D\subset\CC$. Then $f(z,q)$ is iterated per each $q$. One
associates each resulting orbit $\mathcal{O}$ to a connected or
disconnected Julia set, if $\mathcal{O}$ is bounded or not
respectively, so that $q$ may belong or not to $\mathcal{M}$. Thus
$\mathcal{M}$ can be intended as a chart of Julia sets
connectivity, compiled and drawn in the $q-$parameter space.
\begin{figure}[htb]
  \centering
  \fbox{\epsfxsize=3.0cm\epsffile{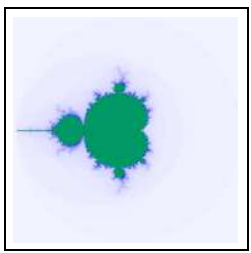}}\caption{The Mandelbrot set for $f(z):z^2+q$.}\label{pic00}
\end{figure}

\begin{question}
What does it happen to $\mathcal{M}$ as $z_0$ is perturbed and it
lies no longer at $c$?
\end{question}
The question here posed was previously introduced in section \S\
6.3 of \cite{Rosa-2005}, where, in the course of related
experiments, we noticed that the connectivity of the Mandelbrot
set $\mathcal{M}$, for the classic quaternionic iterator
$f(h):h^2+q$ $(h,q\in\HH)$, breaks as the seed $h_0$ of each orbit
is distinct from the critical point (in this case at the origin).
Anyway the interpretation of figures with solid 3-D models for the
quaternionic Julia sets looks very hard to evince the question
here posed; therefore we moved to the more simply looking
environment $\CC$, where figures are displayed quite clearer.

\section{MLC conjecture}
The still open MLC (standing for `Mandelbrot Locally Connected')
conjecture states that the Mandelbrot set $\mathcal{M}$ is locally
connected\footnote{If at each point $x\in\mathcal{M}$, every
neighborhood of $x$ includes a connected open neighborhood.},
roughly speaking, the set includes no uni-dimensional filaments.

\section{The journey}
These pages want to enjoy an empirical flavor exclusively, i.e. no
mathematics is achieved and any goal or pretension to deep into
the theory is disclaimed. We will not dare into the proof of it
for sure, because we are not skilled enough for that. This is just
a communication collecting a bunch of figures processed by playing
with the seed $z_0$. Two journeys are shown: given
$R=\mathcal{RE}(q),I=\mathcal{IM}(q)$, we let $q$ range along the
real $R\in[-1.6,0]$ and the imaginary segment $I\in[0,-1.6i]$
respectively, where $\mathcal{M}$ is more visible. We noticed that
as far as one runs away from the origin ($c=0$ for $f(z):z^2+q$)
along both real and imaginary directions, $\mathcal{M}$ loses
progressively connectivity in the peripheral regions, while the
central cardioid is deformed, shrinking again and again. If the
seed is very close to the critical point, disconnectivity is
hardly noticed because only the most marginal periphery of
$\mathcal{M}$ is notched; but, as one moves farther from $c$, the
disconnectivity rate increases, even the cardiod is deformed and
more evident figures are displayed.

\section{One proposal to attack MLC}
We liked to give some (na\"{\i}ve, maybe) impressions: one could
try to attack the MLC conjecture by an inverse approach, that is,
proving that $\mathcal{M}$ is always disconnected when the seed,
for each orbit in the $\mathcal{M}$ experiment for $f(z):z^2+q$ is
not at the origin, according to the original experiment.
\begin{question}
Does the connectivity of $\mathcal{M}$ depend on $z_0$ at the
critical point?
\end{question}
If so, one would try to determine, if possible, a sort of
\emph{disconnectivity rate} decreasing to $0$ along any path
leading to the critical point of the given function $f$. Thus, on
one hand, we pull out the evidence, from the next figures, that
$\mathcal{M}$ is not connected for seeds not lying at the critical
point; on the other hand, from the known figure of $\mathcal{M}$
(see fig. \ref{pic00}), it seems to be locally connected. Perhaps,
things may match together.

\section{Conclusions}
My intention here is just to suggest an idea, an attack to MLC and
know if it might work or not or where it might lead somewhere. In
some cases, like Jacobi's for elliptic integrals, the history of
mathematics showed that an inverted look-up at the problem helped.
Who knows whether the above considerations could actually work?
Without addressing to audience, one would never check it. The next
images are snapshots from a couple of animations we made,
downloadable from our
site\footnote{http://www.malilla.supereva.it/Mirror/Pages/papers.html}
and showing how $\mathcal{M}$ shifts its shape when the seed
$z\neq c$, like in the original experiment.

\bibliographystyle{amsplain}

                        \begin{figure}
                        \centering
                                \subfigure[$R=0.0, I=-0.1$. Infinitesimal disconnections in the
                                peripheral, bottom region; they
                                seem to not happen inside the central cardioid yet.]{\fbox{\epsfxsize=4.5cm\epsffile{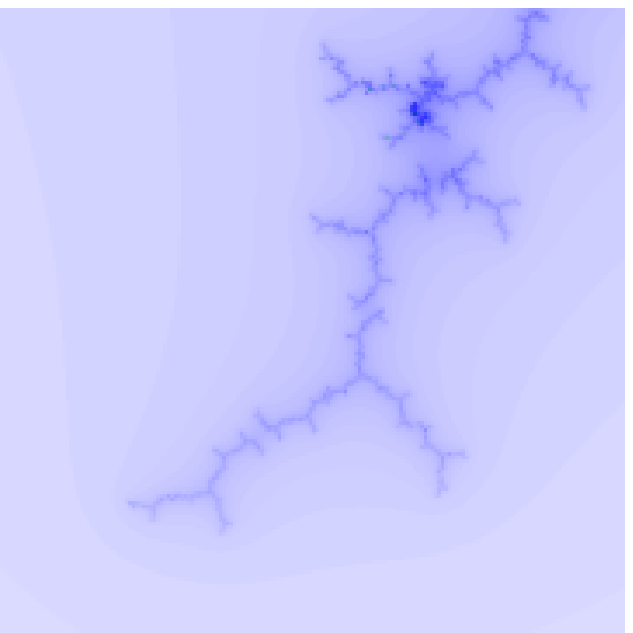}}}\label{pic01}
                                \hspace{0.4cm}
                                \subfigure[$R=0.0, I=-0.2$. Besides the peripheral filaments, one notices a weak deformation in the bigger
                                circle here.]{\fbox{\epsfxsize=4.5cm\epsffile{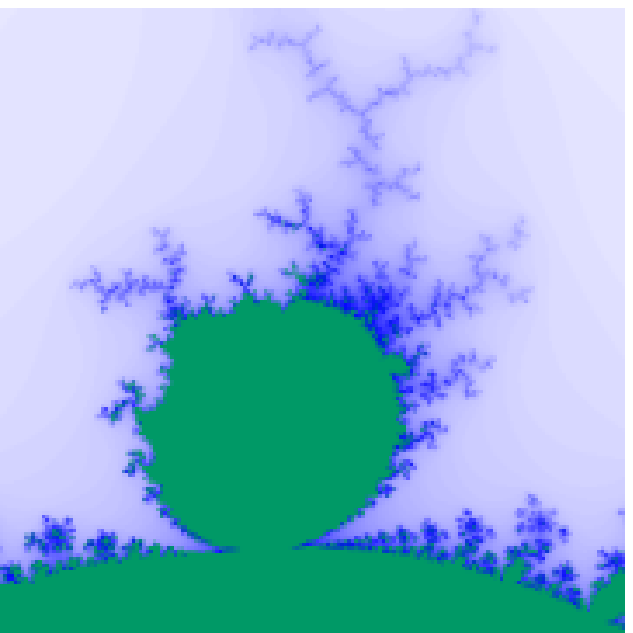}}}\label{pic02}
                                \hspace{0.4cm}
                                \subfigure[$R=0.0, I=-0.4$. Zoom out: the central cardioid is seriously deforming.]
                                {\fbox{\epsfxsize=4.5cm\epsffile{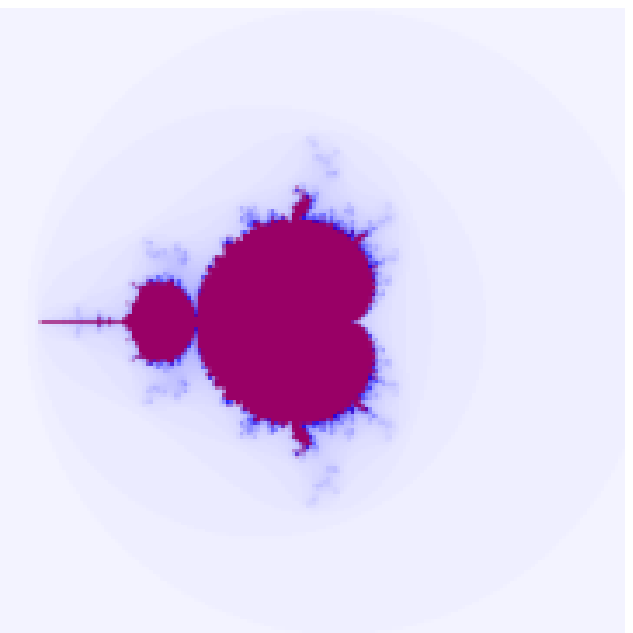}}}\label{pic03}\vspace{1.0cm}\\
                                \subfigure[$R=0.0, I=-0.6$. Further deformation.]
                                {\fbox{\epsfxsize=4.5cm\epsffile{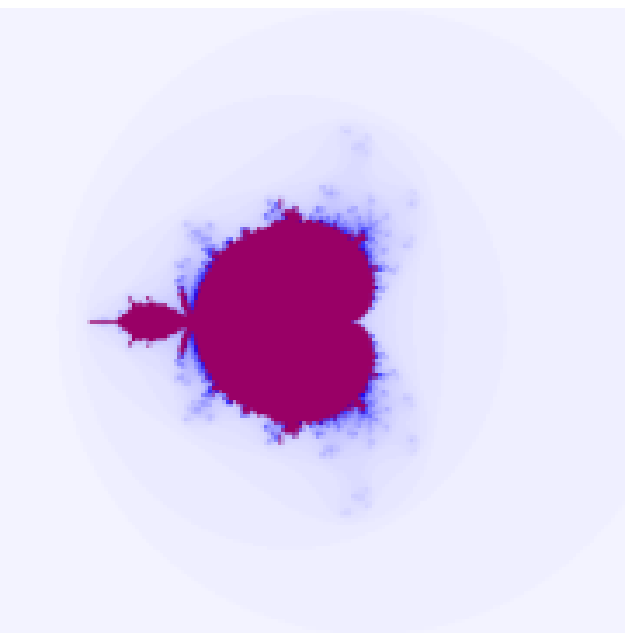}}}\label{pic04}
                                \hspace{0.4cm}
                                \subfigure[The left region of $\mathcal{M}$ shortens into the cardioid.]
                                {\fbox{\epsfxsize=4.5cm\epsffile{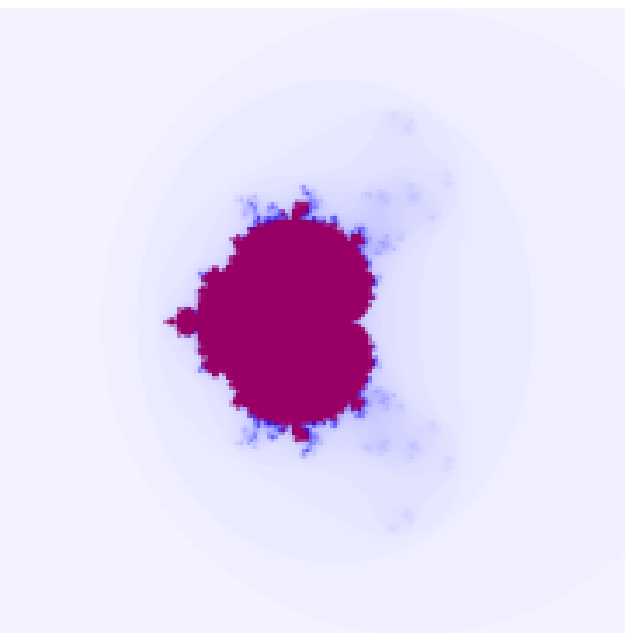}}}\label{pic05}
                                \hspace{0.4cm}
                                \subfigure[$R=0.0, I=-0.9$. The cardioid shape fades away.]
                                {\fbox{\epsfxsize=4.5cm\epsffile{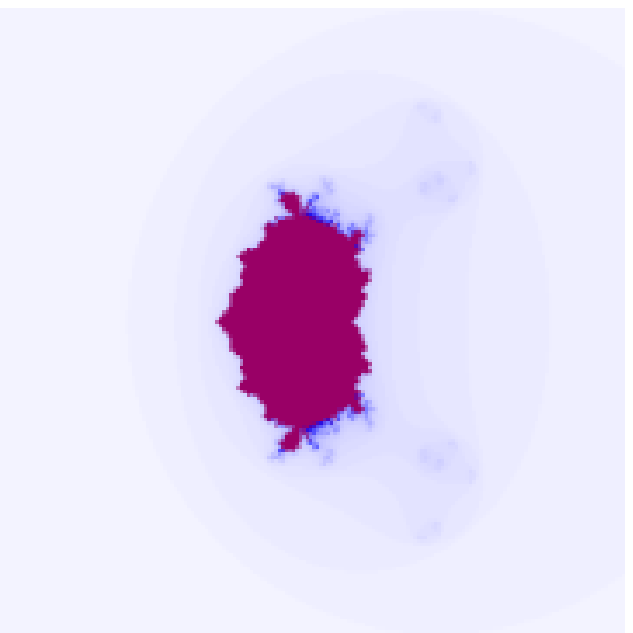}}}\label{pic06}\vspace{1.0cm}\\

                                \subfigure[$R=0.0, I=-1.0$. The previous region
                                turns into two sub-regions whose boundaries intersect at the origin.]
                                {\fbox{\epsfxsize=4.5cm\epsffile{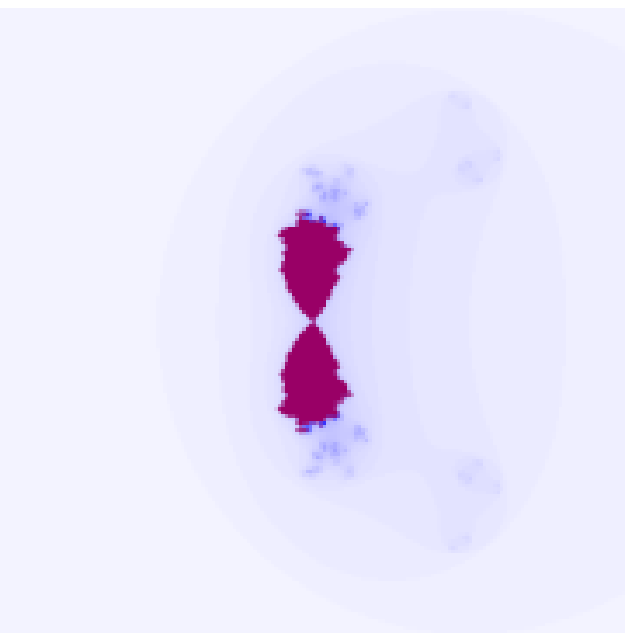}}}\label{pic07}
                                \hspace{0.4cm}
                                \subfigure[A blow up of figure \ref{pic07}.]
                                {\fbox{\epsfxsize=4.5cm\epsffile{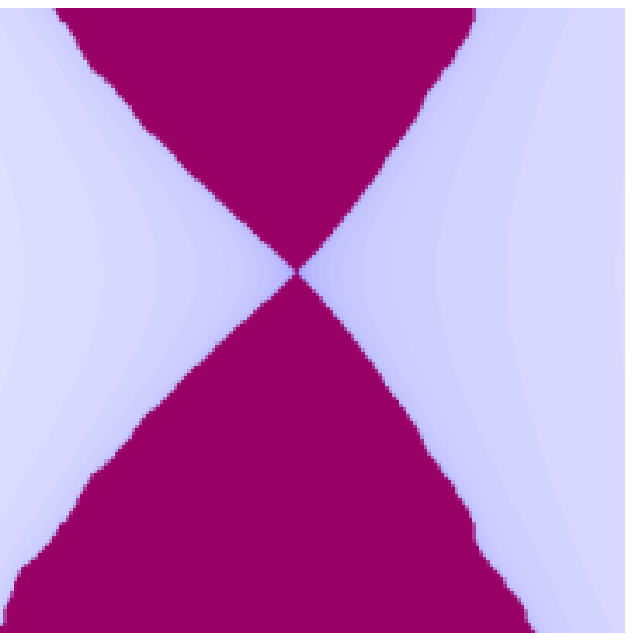}}}\label{pic08}
                                \hspace{0.4cm}
                                \subfigure[$R=0.0, I=-1.05$. The two regions split away.]
                                {\fbox{\epsfxsize=4.5cm\epsffile{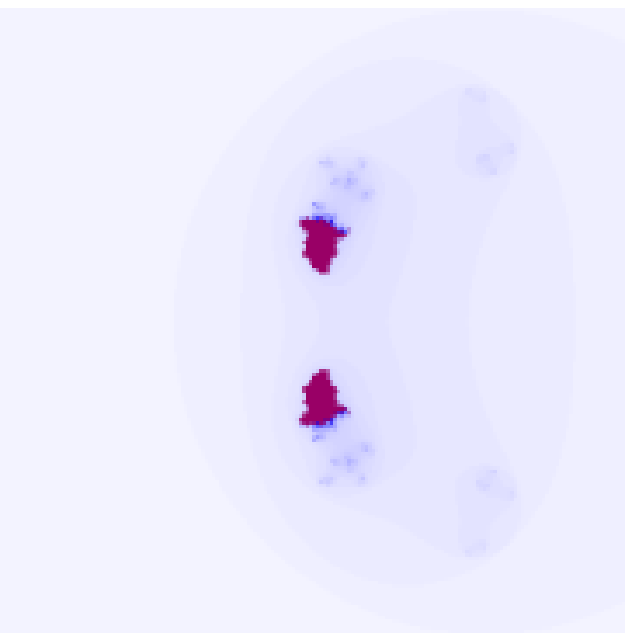}}}\label{pic09}\\
                        \end{figure}

                        \begin{figure}
                        \centering
                                \subfigure[$\mathcal{M}$ appears to be totally disconnected now.]
                                {\fbox{\epsfxsize=4.5cm\epsffile{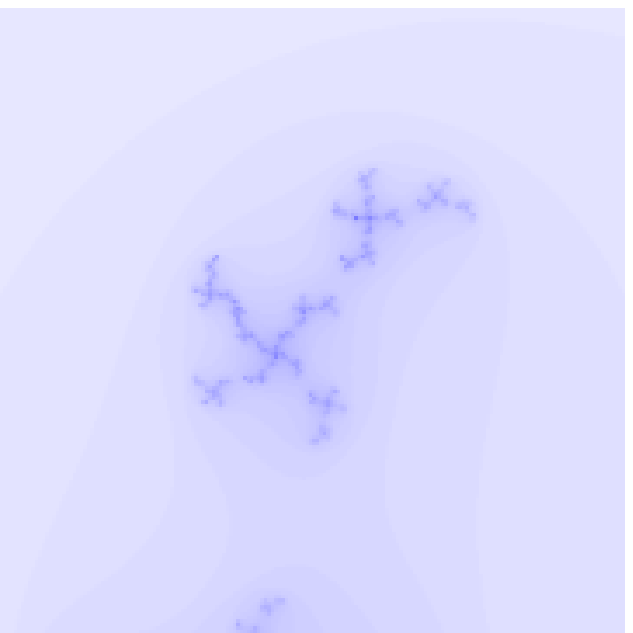}}}\label{pic10}
                                \hspace{0.4cm}
                                \subfigure[$R=0.1, I=0.0$. We start moving along the real segment from 0 to -1.6.]
                                {\fbox{\epsfxsize=4.5cm\epsffile{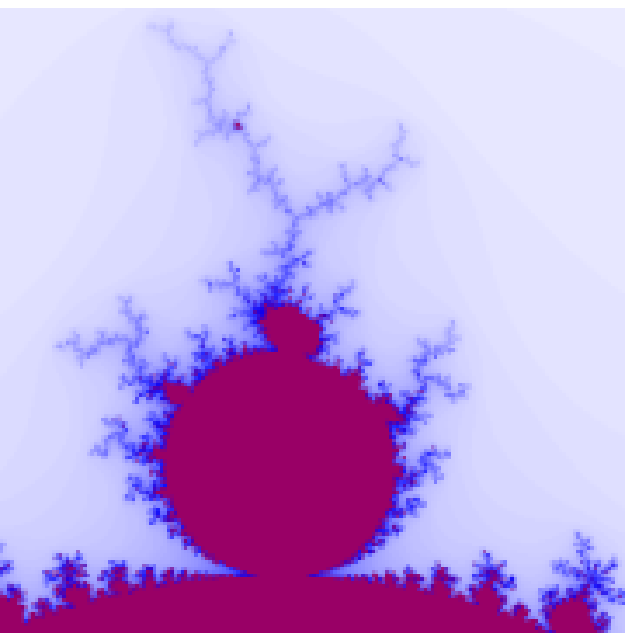}}}\label{pic11}
                                \hspace{0.4cm}
                                \subfigure[$R=0.2, I=0.0$. The same view as in fig. \ref{pic11},but different computation.]
                                {\fbox{\epsfxsize=4.5cm\epsffile{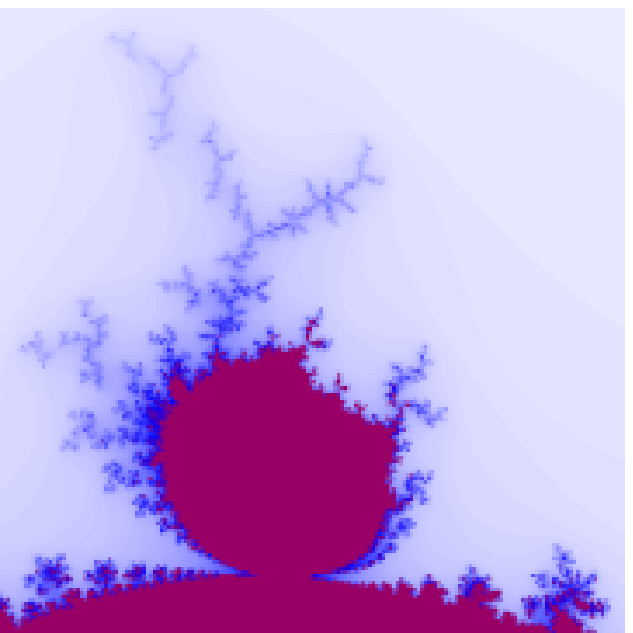}}}\label{pic12}\vspace{1.0cm}\\

                                \subfigure[$\mathcal{M}$ is blowing away again !]
                                {\fbox{\epsfxsize=4.5cm\epsffile{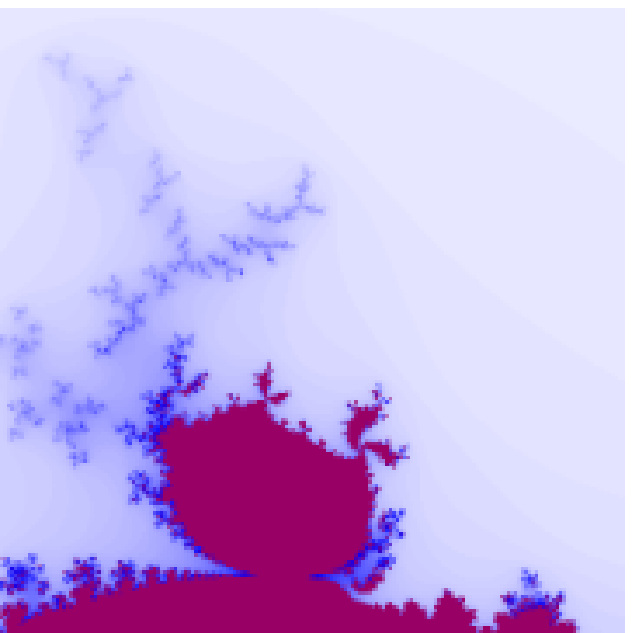}}}\label{pic13}
                                \hspace{0.4cm}
                                \subfigure[For $R=0.6,I=0.0$.]
                                {\fbox{\epsfxsize=4.5cm\epsffile{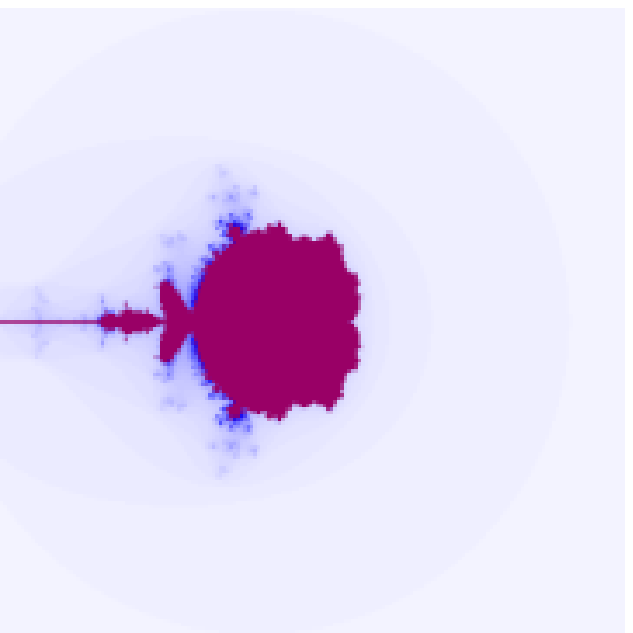}}}\label{pic14}
                                \hspace{0.4cm}
                                \subfigure[$R=1.2, I=0.0$. Speeding deformation up.]
                                {\fbox{\epsfxsize=4.5cm\epsffile{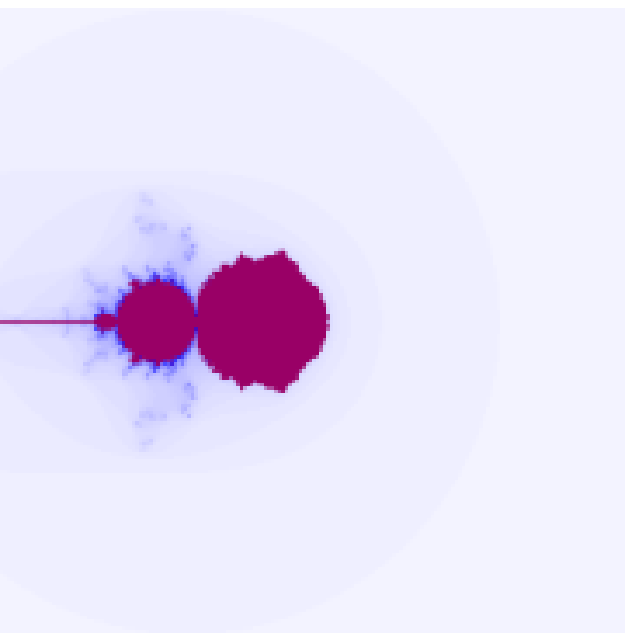}}}\label{pic15}\vspace{1.0cm}\\

                                \subfigure[$R=1.5, I=0.0$. The shrinking.]
                                {\fbox{\epsfxsize=4.5cm\epsffile{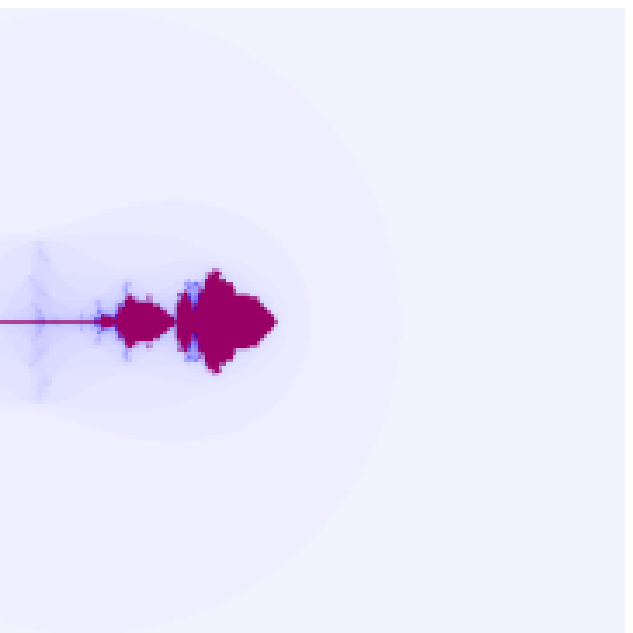}}}\label{pic16}
                                \hspace{0.4cm}
                                \subfigure[$R=1.6, I=0.0$.]
                                {\fbox{\epsfxsize=4.5cm\epsffile{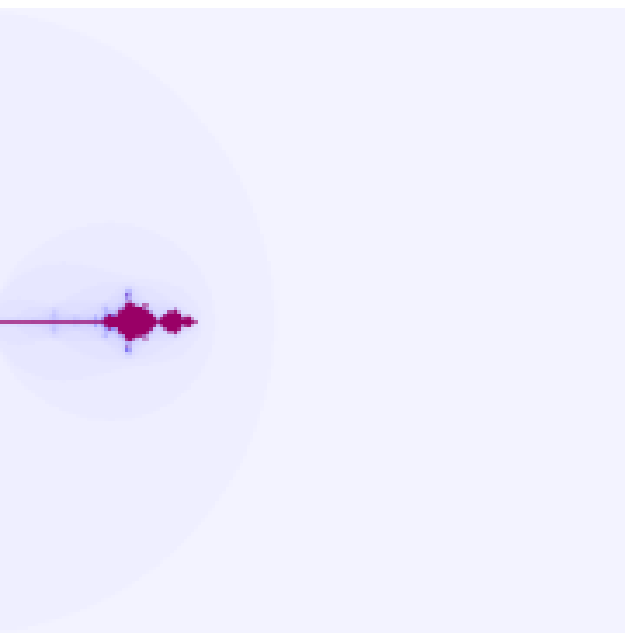}}}\label{pic17}
                                \hspace{0.4cm}
                                \subfigure[$R=1.6, I=0.0$. $\mathcal{M}$ turned into a line.]
                                {\fbox{\epsfxsize=4.5cm\epsffile{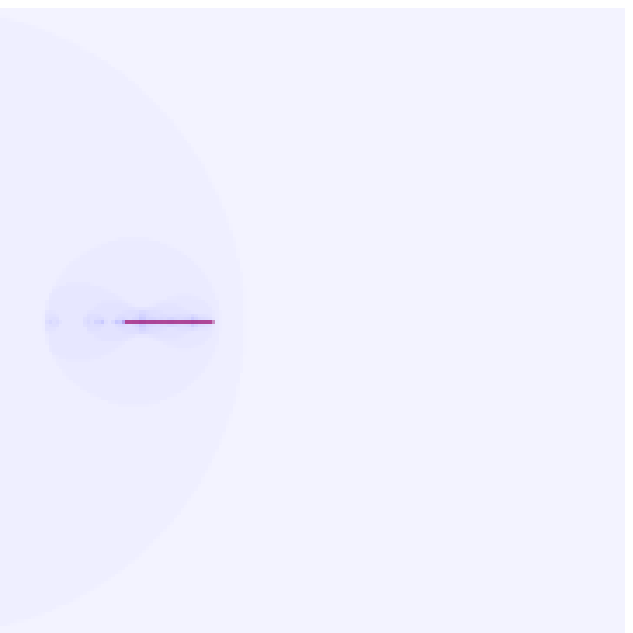}}}\label{pic18}\\
                        \end{figure}

\end{document}